\documentclass[10pt,twoside]{article}
\usepackage{amsmath}
\usepackage{amssymb}
\usepackage{Latex-document}
\usepackage[dvips]{graphicx,color}
\usepackage[dvips]{color}

\newcommand{\dist}{{\rm dist}}

\newcommand{\R}{{\rm Re}}

\newtheorem{theorem}{Theorem}[section]
\newtheorem{proposition}{Proposition}[section]

\newtheorem{statement}{Statement}[section]
\renewcommand{\thesection}{\arabic{section}}
\renewcommand{\thesubsection}{\thesection.\arabic{subsection}.}
\def\addsec{\addtocounter{section}{1} \setcounter{theorem}{0}}

\title{\bf Bifurcations and Strange Attractors\vskip 6mm}

\author{Leonid Shilnikov\vspace*{-0.5cm}\thanks{Institute for Applied Mathematics,
Nizhny Novgorod, Russia. \newline E-mail:shilnikov@focus.nnov.ru}}

\date{\vspace{-8mm}}

\begin{document}

\maketitle

\thispagestyle{first} \setcounter{page}{349}

\begin{abstract}

\vskip 3mm

We reviews the theory of strange attractors and their bifurcations. All known strange attractors may be subdivided
into the following three groups: hyperbolic, pseudo-hyperbolic ones and quasi-attractors. For the first ones the
description of bifurcations that lead to the appearance of Smale-Williams solenoids and Anosov-type attractors is
given. The definition and the description of the attractors of the second class are introduced in the general
case. It is pointed out that the main feature of the attractors of this class is that they contain no stable
orbits. An etanol example of such pseudo-hyperbolic attractors is the Lorenz one. We give the conditions of their
existence. In addition we present a new type of the spiral attractor that requires countably many topological
invariants for the complete description of its structure. The common property of quasi-attractors and
pseudo-hyperbolic ones is that both admit homoclinic tangencies of the trajectories. The difference between them
is due to quasi-attractors may also contain a countable subset of stable periodic orbits. The quasi-attractors are
the most frequently observed limit sets in the problems of nonlinear dynamics. However, one has to be aware that
the complete qualitative analysis of dynamical models with homoclinic tangencies cannot be accomplished.

\vskip 4.5mm

\end{abstract}

\vskip 12mm

\section*{1. Introduction}\label{Section1} \addsec

\vskip-5mm \hspace{5mm}

One of the starling discoveries in the XX century was the discovery of dynamical chaos. The finding added a
fascinatingly novel type of motions
--- chaotic oscillations to the catalogue of the accustomed, in nonlinear
dynamics, ones such as steady states, self-oscillations and
modulations. Since then  many problems of contemporary exact and engineering
sciences that are modelled within the framework
of differential equations have obtained the adequate mathematical
description. Meanwhile this also set up a question: the trajectory of what
kind are to be connected to dynamical chaos in systems with 3D and higher
dimensions of the phase space.

The general classification of the orbits in dynamical systems is
due to Poincar\'e and Birkhoff. As far as chaos is concerned in particular
a researcher is interested most of all in non-wandering trajectories
which are, furthermore, to be self-liming.
In the hierarchy of all orbits of dynamical system the high goes to
the set of the unclosed trajectories so-called stable in {\it Poisson's sense}. The
feature of them is that the sequence composed of the differences of the
subsequent Poincar\'e recurrence times of a trajectory (as it gets back
into the vicinity of the initial point) may have no upper bound.
In the case where the return times of the trajectory are bounded
Birkhoff suggested to call such it a {\em recurrent} one.
A partial subclass of recurrent
trajectories consists of almost-periodic ones which, by definition, possess
a set of almost-periods. Since the spectrum of the set is discrete it
follows that the dynamics of a system with the almost-periodic trajectories
must be simple. Moreover, the almost-periodic trajectories are also quasi-periodic,
the latter ones are associated with the regime of modulation. The closure
of a recurrent trajectory is called a {\em minimal set}, and that of
a Poisson stable trajectory is called a {\em quasi-minimal} one. The
last one contains also the continuum of Poisson stable trajectories which
are dense in it. The quasi-minimal set may, in addition, contain some closed
invariant subsets such as equilibria, periodic orbits, ergodic
invariant tori with quasi-periodic covering and so on. All this is the
reason why  Poincar\'e recurrent times of a Poisson stable
trajectory are unbounded: it may linger in the neighborhoods of
the above subsets for rather long times before it passes by the initial point.
In virtue of such unpredictable behavior it seems quite reasonable the
the the role of dynamical chaos orbits should be assigned to the Poisson stable
trajectories.

Andronov was the first who had raised the question about the
correspondence between the classes of trajectories of dynamical
systems and the types of observable motions in nonlinear dynamics
in his work ``Mathematical problems of auto-oscillations".
Because he was motivated to explain the nature of
self-oscillations he repudiated the Poisson stable trajectory
forthwith due to their irregular behavior. He expressed the
hypothesis that a recurrent trajectory stable in Lyapunov sense
would be almost-periodic; he also proposed to Markov to confirm
it. Markov proved a stronger result indeed; namely that a
trajectory stable both in Poisson and Lyapunov (uniformly) senses
would be an almost-periodic one. This means that the Poisson
stable trajectories must be unstable in Lyapunov sense to get
associated to dynamical chaos. After that it becomes clear that
all trajectories in a quasi-minimal set are to be of the saddle
type. The importance of such quasi-minimal sets for nonlinear
dynamics as the real objects was inferred in the explicit way by
Lorenz in 1963. He presented the set of equations, known today as
the Lorenz model, that possessed an attracting quasi-minimal set
with the unstable trajectory behavior. Later on, such sets got
named {\em strange attractors} after Ruelle and Takens.

Here we arrive at the following problem: how can one establish the existence
of the Poisson stable trajectories in the phase space of a system? Furthermore,
the applicability of the system as a nonlinear dynamics model
requires that such trajectories shall persist under small smooth perturbations.
Undoubtedly the second problem is a way complex than that of finding periodic
orbits. Below we will yield a list of conditions that guarantee the existence
of unstable trajectories stable in Poisson sense.

The most universal criterion of the existence of Poisson stable trajectories is the
presence, in a system, of a hyperbolic saddle periodic orbit whose stable and
unstable manifolds cross transversally along the homoclinic orbit. This structure
implies that the set $N$ of all of the orbits remaining in its small
neighborhood consists of only unstable ones. Moreover, the periodic
orbits are dense in $N$, so are the trajectories homoclinic to them,
besides the continuum of unclosed trajectories stable in Poisson sense.
Generally speaking wherever one can describe the behavior of the trajectories in
terms of symbolic dynamics (using either the Bernoulli sub-shift or the
Markov topological chains), the proof of the existence of the Poisson
stable trajectories becomes trivial. However, the selected hyperbolic
sets by themselves are unstable. Nevertheless, their presence in the phase
space means the complexity of the trajectory behavior even though they
are no part of the strange attractor.

Early sixtieth were characterized the rapid development of the
theory of structural stability initiated in the works of Anosov
and Smale. Anosov was able to single out the class of the systems
for which the hyperbolicity conditions hold in the whole phase
space. Such flows and cascade have been named the Anosov systems.
Some examples of the Anosov systems include geodesic flows on
compact smooth manifolds of a negative curvature \cite{AN67}. It
is well-known that such flow is conservative and its set of
non-wandering trajectories coincides with the phase space.
Another  example of the Anosov diffeomorphism is a mapping of an
$n$-dimensional torus
\begin{equation}\label{p4.7}
\bar \theta = A \theta + f(\theta), \qquad  \mbox {mod    } 1,
\end{equation}
where $A$ is a matrix with integer entries other than $1$, $\det | A | =1$,
the eigenvalues of $A$ do not lie on the unit circle,
and  $f(\theta)$ is a periodic function of period $1$.

The condition of hyperbolicity of (\ref{p4.7}) may be easily verified for
one pure class of diffeomorphisms of the kind:
\begin{equation}\label{p4.9}
\bar \theta = A \theta, \qquad  \mbox {mod    } 1,
\end{equation}
which are the algebraic hyperbolic automorphisms of the torus.
Automorphism~(\ref{p4.9}) are conservative systems whose set $\Omega$ of
non-wandering trajectories coincides with the torus ${\Bbb T}^{n}$ itself.

Conditions of structural stability of high-dimensional systems was
formulated by Smale \cite{SM67}. These conditions are in the
following: A system must satisfy both {\it Axiom~A} and
the strong transversality condition.

Axiom~A requires that:
\begin{enumerate}
\item[1A]
the non-wandering set $\Omega$ be hyperbolic;
\item[1B]
$\Omega = \overline{ \mbox{Per}}$. Here $\mbox{Per}$ denotes the set of
periodic points.
\end{enumerate}
Under the assumption of Axiom~A the set $\Omega$ can be represented by
a finite union of non-intersecting, closed, invariant, transitive sets
$\Omega_{1}, \ldots \Omega_{p}$. In the case of cascades, any such
$\Omega_{i}$ can be represented by a finite number of sets having these
properties which are mapped to each other under the action of the
diffeomorphism. The sets  $\Omega_{1}, \ldots \Omega_{p}$ are called
{\it basis sets}.

The basis sets of Smale systems (satisfying the enumerated conditions) may
be of the following three types: attractors, repellers and saddles.
Repellers are the basis sets which becomes attractors in backward time.
Saddle basis sets are such that may both attract and repel outside
trajectories. A most studied saddle basis sets are one-dimensional in the case
of flows
and null-dimensional in the case of cascades. The former ones are
homeomorphic to the suspension over topological Markov chains; the latter
ones are homeomorphic to simple topological Markov chains [Bowen \cite{BO70}].

Attractors of Smale systems are called {\it hyperbolic}.
The trajectories passing sufficiently close
to an attractor of a Smale system, satisfies the condition

$$
\dist ( \varphi(t,x), A) < k e^{-\lambda t}, \qquad t\ge 0
$$
where $k$ and $\lambda$ are some positive constants. As we have said earlier
these attractors are transitive. Periodic, homo- and heteroclinic trajectories
as well as Poisson-stable ones are everywhere dense in them. In particular,
we can tell one more of their peculiarity: the unstable manifolds of all points
of such an attractor lie within it, i.e., $W^{s}_{x} \in A$ where $x \in A$.
Hyperbolic attractors may be smooth or non-smooth manifolds, have a
fractal structure, not locally homeomorphic to a direct product of a disk and
a Cantor set and so on.

Below we will discuss a few hyperbolic attractors which might be curious
for nonlinear dynamics. The first example of such a hyperbolic attractor
may be the Anosov torus ${\Bbb T}^{n}$ with a hyperbolic structure on it.
The next example of hyperbolic attractor was constructed by Smale on a
two-dimensional torus by means of a ``surgery'' operation over the
automorphism of this torus with a hyperbolic structure. This is the so-called
$DA$-({\it derived from Anosov}) diffeomorphism. Note that
the construction of such attractors is designed as that of minimal sets
known from the Poincar\'e-Donjoy theory in the case of ${\Bbb C}^{1}$-smooth
vector fields on a two-dimensional torus.

Let us consider a solid torus $\Pi \in {\Bbb R}^{n}$, i.e., ${\Bbb
T}^{2}= {\Bbb D}^{2}\times {\Bbb S}^{1}$ where ${\Bbb D}^{2}$ is a
disk and ${\Bbb S}^{1}$ is a circumference. We now expand ${\Bbb
T}^{2}$ m-times (m is an integer) along the cyclic coordinate on
${\Bbb S}^{1}$ and shrink it $q$-times along the diameter of
${\Bbb D}^{2}$ where $q \le 1/m$. We then embed this deformated
torus $\Pi_{1}$ into the original one so that its intersection
with ${\Bbb D}^{2}$ consists of $m$-smaller disks. Repeat this
routine with $\Pi_{1}$ and so on. The set $\Sigma =
\in_{i=1}^{\infty} \Pi_{i}$ so obtained is called a
 Witorius-Van Danzig solenoid.
Its local structure may be represented as the direct product of an interval
and a Cantor set.
Smale also observed that Witorius-Van Danzig solenoids may have
a hyperbolic structures, i.e.,  be hyperbolic attractors of
diffeomorphisms on solid tori. Moreover, similar attractors can be realized
as a limit of the inverse spectrum of the expanding cycle map \cite{W67}
$$
\theta = m \theta, \qquad \mbox {mod    }1.
$$
The  peculiarity of such solenoids is that they are expanding solenoids.
Generally speaking, an expanding solenoid is called a hyperbolic attractor
such that its dimension coincides with the dimension of the unstable manifolds
of the points of the attractor. Expanding solenoids were studied by Williams
\cite{W74}
who showed that they are generalized (extended) solenoids. The construction of
generalized solenoids is similar to that of minimal sets of
limit-quasi-periodic trajectories. Note that in the theory of sets of
limit-quasi-periodic functions the  Wictorius-Van Danzig solenoids   are
quasi-minimal sets. Hyperbolic solenoids are called the Smale-Williams
solenoids.

We remark also on an example of a hyperbolic attractor of a
diffeomorphism on a two-dimensional sphere, and, consequently, on
the plane, which was built by Plykin \cite{PSS91}. In fact, this
is a diffeomorphism of a two-dimensional torus projected onto a
two-dimensional sphere. Such a diffeomorphism, in the simplest
case, possesses four fixed points, moreover all of them are
repelling.

\vspace*{-2mm}

\section*{2. Birth of hyperbolic attractors} \addsec

\vskip-5mm \hspace{5mm}

Let us now pause to discuss the principal aspects related to the transition from Morse-Smale systems to systems
with hyperbolic attractors. The key moment here is a global bifurcation of disappearance of a periodic trajectory.
Let us discuss this bifurcation in detail following the paper of Shilnikov and Turaev \cite{BLUE}.

Consider a ${\Bbb C}^{r}$-smooth one-parameter family of dynamical systems
$X_{\mu}$ in ${\Bbb R}^{n}$. Suppose that the flow has a periodic orbit
$L_{0}$ of the saddle-node type at $\mu =0$. Choose a neighbourhood
$U_{0}$ of $L_{0}$ which is a solid torus partitioned by the $(n-1)$-dimensional
strongly stable manifold $W^{ss}_{L_{0}}$ into two regions: the node region
$U^{+}$ where all trajectories tend to $L_{0}$ as $t \to +\infty$, and the
saddle region $U^{-}$  where the two-dimensional unstable manifold
$W^{u}_{L_{0}}$ bounded by $L_{0}$ lies.
Suppose that all of the trajectories of $W^{u}_{L_{0}}$ return to
$L_{0}$ from the node region $U^{+}$ as $t \to +\infty$ and do not
lie in $W^{ss}$. Moreover, since any trajectory of $W^{u}$ is
bi-asymptotic to $L_{0}$, $\bar W^{u}_{L_{0}}$ is compact.

Observe that systems close to $X_{0}$ and having a simple saddle-node periodic
trajectories close $L_{0}$ form a surface $B$ of codimension-1 in the space of
dynamical systems. We assume also that the family $X_{\mu}$ is transverse to
$B$. Thus, when $\mu < 0$, the orbit $L_{0}$ is split into two periodic orbits,
namely: $L^{-}_{\mu}$ of the saddle type and stable $L^{+}_{\mu}$.
When  $\mu > 0$  $L_{0}$ disappears.

It is clear that $X_{\mu}$ is a Morse-Smale system in a small neighbourhood
$U$ of the set $W^{u}$ for all small $\mu <0$. The non-wandering set here
consists of the two periodic orbits $L^{+}_{\mu}$ and $L^{-}_{\mu}$. All
trajectories of $U \backslash W^{s}_{L^{-}_{\mu}}$ tend to $L_{\mu}^{+}$
as $t \to  +\infty$. At $\mu=0$ all trajectories on $U$ tend to $L_{0}$.
The situation is more complex when $\mu >0$.

The Poincar\'e map to which the problem under consideration is reduced, may be
written in the form

\begin{equation}\label{p6.1}
\begin{array}{l}
\bar x = f(x, \theta, \mu), \\
\bar \theta = m \theta + g(\theta) + \omega + h(x, \theta, \mu),
 \qquad \mbox {mod  }1,
\end{array}
\end{equation}
where $f,\, g$ and $h$ are periodic functions of $\theta$. Moreover,
$\| f \|_{{\Bbb C}^{1}} \to 0$ and $\|h\|_{{\Bbb C}^{1}} \to 0$
as $\mu \to 0$, $m$ is an integer and $\omega$ is a parameter defined in the
set $[0,1 )$.  Diffeomorphism~(\ref{p6.1}) is defined in a solid-torus
${\Bbb D}^{n-2} \times {\Bbb S}^{1}$, where ${\Bbb D}^{n-2}$ is a disk
$\|x \| <r$, $r>0$
Observe that (\ref{p6.1}) is a strong contraction along $x$.
Therefore, mapping~(\ref{p6.1}) is close to the degenerate map

\begin{equation}\label{p6.2}
\begin{array}{l}
\bar x = 0, \\
\bar \theta = m \theta + g(\theta) + \omega, \qquad \mbox {mod  }1.
\end{array}
\end{equation}
This implies that its dynamics is determined by the circle map

\begin{equation}\label{p6.3}
\bar \theta = m \theta + g(\theta) + \omega, \qquad \mbox {mod  }1,
\end{equation}
where $0 \le \omega < 1$. Note that in the case                   of the
flow in ${\Bbb R}^{3}$, the integer $m$ may assume the values $0,\, 1$.

\begin{theorem}\label{t6.1}
If $m=0$ and if
$$
\max \| g^{\prime}(\theta)\| < 1,
$$
then for sufficiently small $\mu > 0$, the original flow has a periodic  orbit
both of which length and period tend to infinity as $\mu \to 0$.
\end{theorem}

This is the ``blue sky catastrophe''. In the case  where
$m=1$, the closure $\bar W^{u}_{L_{0}}$ is a two-dimensional torus. Moreover,
it is smooth provided that (\ref{p6.1}) is a diffeomorphism. In the case where $m=-1$
$\bar W^{u}_{L_{0}}$ is a Klein bottle, also smooth if (\ref{p6.1}) is a
diffeomorphism. In the case of the last theorem $\bar W^{u}_{L_{0}}$ is  not a
manifold.

In the case of ${\Bbb R}^{n}$ $(n\ge 4)$ the constant $m$ may be any integer.

\begin{theorem}
Let $|m | \ge 2$ and let $| m + g^{\prime}(\theta) > 1$. Then for all
$\mu > 0$ sufficiently small, the Poincar\'e map~(\ref{p6.1}) has a
hyperbolic attractor homeomorphic to the Smale-Williams solenoid, while
the original family has a hyperbolic attractor homeomorphic to a
suspension over the Smale-Williams solenoid.
\end{theorem}

The idea of the use of the saddle-node bifurcation to produce  hyperbolic
attractors may be extended onto that of employing the bifurcations of an
invariant torus.  We are not developing here the theory of such bifurcations but
restrict ourself by consideration of a modelling situation.

Consider a one-parameter family of smooth dynamical systems
$$
\dot x = X(x, \mu)
$$
which possesses an invariant $m$-dimensional torus ${\Bbb T}^{m}$ with a quasi-periodic
trajectory at $\mu =0$.  Assume that the vector field may be recast as
\begin{equation}\label{p6.4}
\begin{array}{l}
\dot y = C(\mu) y, \\
\dot z = \mu + z^{2}, \\
\dot \theta = \Omega(\mu)
\end{array}
\end{equation}
in a neighborhood of ${\Bbb T}^{m}$. Here, $z \in {\Bbb R}^{1}$,
$y \in {\Bbb R}^{n-m-1}$, $\theta \in {\Bbb T}^{m}$ and
$\Omega(0)=(\Omega_{1}, \ldots, \Omega_{m})$. The matrix $C(\mu)$ is
stable, i.g., its eigenvalues lie to the left of the imaginary axis in the
complex plane. At $\mu =0$ the equation of the torus is $y=0$, the equation
of the unstable  manifold $W^{u}$ is $y=0, z >0$, and that of the strongly
unstable manifold $W^{ss}$ partitioning the neighborhood of ${\Bbb T}^{m}$
into a node and a saddle region, is $z=0$. We assume also that all of the
trajectories of the unstable manifold $W^{u}$ of the torus come back to it
as $t \to +\infty$. Moreover they do not lie in $W^{ss}$. On a cross-sections
transverse to $z=0$ the associated Poincar\'e map may be written in the
form

\begin{equation}\label{p6.5}
\begin{array}{l}
\bar y = f(y, \theta, \mu), \\
\bar \theta = A \theta + g(\theta) + \omega + h(x, \theta, \mu),
 \qquad \mbox {mod  }1,
\end{array}
\end{equation}
where A is an integer matrix, $f,\, g,$ and $h$ are $1$-periodic
functions of $\theta$. Moreover $\|f\|_{{\Bbb C}^{1}} \to 0$ and
$\|h\|_{{\Bbb C}^{1}} \to 0$ as $\mu \to 0$, $\omega=(\omega_{1},
\cdots, \omega_{m})$ where $0\le \omega_{k} < 1$.

Observe that the restriction of the Poincar\'e map on the invariant torus is
close to the shortened map
\begin{equation}\label{p6.6}
\bar \theta = A \theta +g(\theta) + \omega, \qquad \mbox {mod  }1.
\end{equation}
This implies, in particular, that if (\ref{p6.6}) is an Anosov map for
all $\omega$ (for example when the eigenvalues of the matrix $A$ do not
lie on the unit circle of the complex plane, and $g(\theta)$ is small),
then the restriction of the Poincar\'e  map  is also an Anosov map for
all $\mu >0$.
Hence, we arrive at the following statement

\begin{proposition}
If the shortened map is an Anosov map for all small $\omega$, then for all
$\mu > 0$ sufficiently small, the original flow possesses a hyperbolic
attractor which is topologically conjugate to the suspension over the
Anosov diffeomorphism.
\end{proposition}

The birth of hyperbolic attractors may be proven not only in the case where
the shortened map is a diffeomorphism. Namely, this result holds true if the
shortened map is expanding. A map
is called expanding of the length if any tangent vector field grows
exponentially under the action of the differential of the map. An example is
the algebraic map
$$
\bar \theta = A \theta, \qquad \mbox {mod   } 1,
$$
such that the spectrum of the integer matrix $A$ lies strictly outside the
unit circle, and any neighboring map is also expanding. If
$\| (G^{\prime}(\theta))^{-1} \| < 1$, where $G=A+g(\theta)$, it follows then that
the shortened map
\begin{equation}\label{p6.7}
\bar \theta =  = \omega  + A \theta +g(\theta),
\qquad \mbox {mod  }1,
\end{equation}
is an expansion for all $\mu >0$.

Shub \cite{SHUB78} established that expanding maps are structurally stable.
The study of expanding maps and their connection to smooth diffeomorphisms
was continued by Williams [Williams \cite{W70}]. Using the result of his
work we come to the following result which is analogous to our
theorem~\ref{t6.1}, namely

\begin{proposition}
If $\| (G^{\prime}(\theta))^{-1} \| < 1$, then for all small $\mu >0$, the
Poincar\'e map possesses a hyperbolic attractor locally homeomorphic to a
direct product of ${\Bbb R}^{m+1}$ and a Cantor set.
\end{proposition}

An endomorphism of a torus is called {\it an Anosov covering} if there exists
a continuous decomposition of the tangent space into the direct sum of stable
and unstable submanifolds just like in the case of the Anosov map (the
difference is that the Anosov covering is not a one-to-one map, therefore,
it is not a diffeomorphism). The map~(\ref{p6.6}) is an Anosov covering
if, we assume, $| \det A | >1$ and if $g(\theta)$ is sufficiently small. Thus,
the following result is similar to the previous proposition

\begin{proposition}
If the shortened map~(\ref{p6.6}) is an Anosov covering for all $\omega$, then
for all small $\mu >0$ the original Poincar\'e map possesses a hyperbolic
attractor locally homeomorphic to a direct product of ${\Bbb R}^{m+1}$ and a
Cantor set.
\end{proposition}

In connection with the above discussion we can ask what other hyperbolic
attractors may be generated from Morse-Smale systems?

Of course there are other scenarios of the transition from a Morse-Smale
system to a system with complex dynamics, for example, through
$\Omega$-explosion, period-doubling cascade, etc. But these
bifurcations do not lead explicitly to the appearance of hyperbolic
strange attractors.

\vspace*{-2mm}

\section*{3. Lorenz attractors} \addsec

\vskip-5mm \hspace{5mm}

In 1963 Lorenz \cite{L63} suggested the model:
\begin{eqnarray}\label{p7.1}
\begin{array}{ll}
{\dot x} = -\sigma(x - y) ,\\
{\dot y} = r x - y - xz ,\\
{\dot z} = -b z +x y \\
\end{array}
\end{eqnarray}
in which he discovered numerically a vividly chaotic behaviour of the
trajectories when $\sigma=10$, $b=8/3$ and $r=28$.
the important conclusion has been derived from the mathematical studies
of the Lorenz model:   simple models of nonlinear dynamics may have
strange attractors.

Like hyperbolic attractors, periodic as well as homoclinic
orbits are everywhere dense in the Lorenz attractor. Unlike hyperbolic
attractors the Lorenz
attractor is structurally unstable. This is due to the embedding of a
saddle equilibrium state with a one-dimensional unstable manifold
into the attractor.  Nevertheless, under small smooth perturbations
stable periodic orbits do not arise. Moreover, it became obvious that
such strange attractors may be obtained through a finite number of
bifurcations. In particular, in the Lorenz model (due to its specific feature:
it has the symmetry group $(x,y,z) \leftrightarrow (-x,-y,z)$)
such a route consists of three steps only.

Below we present a few statements concerning the description of the structure
of the Lorenz attractor as it was done in \cite{ABSH77,ABSH83}. The fact that we are
considering only three-dimensional systems is not important, in principle,
because the general case where only one characteristic value is positive for the
saddle while the others have negative real parts, and the value least with the
modulus  is real, the result is completely similar to the three-dimensional case.
Let $B$
denote the Banach space of ${\Bbb C}^{r}$-smooth dynamical systems ($r \ge 1$)
with the ${\Bbb C}^{r}$-topology, which are specified on  a smooth
three-dimensional manifold $M$. Suppose that in the domain $U \subset B$ each
system $X$ has an equilibrium state $O$ of the saddle type. In this case
the inequalities $\lambda_{1} < \lambda_{2} < 0 < \lambda_{3}$ hold for the
roots $\lambda_{i} =\lambda_{i}(X)$, $i=1,2,3$ of the characteristic equation
at $O$, and the saddle value $\sigma(X)=\lambda_{2}+ \lambda_{3} > 0$. A
stable two-dimensional manifold of the saddle will be denoted by
$W^{s}=W^{s}(X)$ and the unstable one, consisting of $O$ and two trajectories
$\Gamma_{1,2} =\Gamma_{1,2}(X)$ originating from it, by $W^{u}=W^{u}(X)$. It is
known that both $W^{s}$ and $W^{u}$ depend smoothly on $X$  on each compact
subset. Here it is assumed that in a certain local map
$V=\{(x_{1}, x_{2}, x_{3})\}$, containing $O$, $X$ can be  written in the
form
\begin{equation}\label{p7.3}
\dot x_{i} = \lambda_{i} x_{i} + P_{i}(x_{1}, x_{2}, x_{3}),
\qquad i=1,2,3.
\end{equation}
Suppose that the following conditions are satisfied for the system
$X_{0} \subset U$ (see Fig.1):
\begin{figure}
\begin{center}
\includegraphics[width=3 in,height=2.5 in]{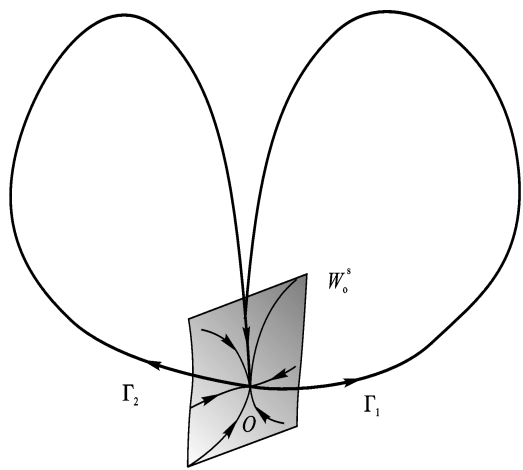}
\end{center}
\caption{Homoclinic butterfly} \label{fig1}
\end{figure}
\begin{enumerate}
\item
$\Gamma_{i}(X_{0}) \subset W^{s}(X_{0})$, $i=1,2$ , i.e., $\Gamma_{i}(X_{0})$ is
doubly asymptotic to $O$.
\item
$\Gamma_{1}(X_{0})$ and $\Gamma_{2}(X_{0})$ approach to $O$ tangentially to
each other.
\end{enumerate}

The condition $\lambda_{1} < \lambda_{2}$ implies that the
non-leading manifold $W^{ss}_{O}$ of $W^{s}_{O}$, consisting of
$O$ and the two trajectories tangential to the axis $x_{1}$ at
the point $O$, divides $W^{s}_{O}$ into two open domains:
$W^{s}_{+}$ and $W^{s}_{-}$. Without loss of generality we may
assume that $\Gamma_{i}(X_{0}) \subset W^{s}_{+}(X_{0})$, and
hence $\Gamma_{i}$ is tangent to the positive semiaxis $x_{2}$.
Let $v_{1}$ and $v_{2}$ be sufficiently small neighborhoods of
the separatrix ``butterfly'' $\bar \Gamma = \overline{\Gamma_{1}
\cup O \cup \Gamma_{2}}$. Let ${\cal M_{i}}$ stand for the
connection component of the intersection of
$\overline{W^{s}_{+}(X_{0})}$ with $v_{i}$, which contains
$\Gamma_{i}(X_{0})$. In the general case ${\cal M}_{i}$ is a
two-dimensional ${\Bbb C}^{0}$-smooth manifold homeomorphic ether
to a cylinder or to a M\"obius band. The general condition lies in
the fact that certain values $A_{1}(X_{0})$ and $A_{2}(X_{0})$,
called the separatrix values, should not be equal to zero.

It follows from the above assumptions that $X_{0}$ belongs to the bifurcation
set $B^{2}_{1}$ of codimension two, and $B^{2}_{1}$ is the intersection of
two bifurcation surfaces $B^{1}_{1}$ and $B^{1}_{2}$ each of codimension one,
where $B^{1}_{i}$ corresponds to the separatrix loop
$\bar \Gamma_{i} = \overline{O \cup \Gamma_{i}}$. In such a situation it is
natural to consider a two-parameter family of dynamical systems $X(\mu)$,
$\mu = (\mu_{1}, \mu_{2})$, $|\mu| <\mu_{0}$, $X(0)=X_{0}$, such that
$X(\mu)$ intersects with $B^{2}_{1}$ only along $X_{0}$ and only for $\mu=0$.
It is also convenient to assume that the family $X(\mu)$ is transverse to
$B^{2}_{1}$. By transversality we mean that for the system $X(\mu)$ the
loop  $\Gamma_{1}(X(\mu))$ ``deviates'' from $W^{s}_{+}(X(\mu))$ by a value
of the order of $\mu_{1}$, and
the loop  $\Gamma_{2}(X(\mu))$ ``deviates'' from $W^{s}_{+}(X(\mu))$ by a value
of the order of $\mu_{2}$.

It is known from \cite{SH70} that the above assumptions imply that
in the transition to a system close to $X_{0}$ the separatrix loop
can generate only one periodic orbit which is of the saddle type.
Let us assume, for certainty, that the loop $\Gamma_{1}(X_{0})
\cup O$ generates a periodic orbit $L_{1}$ for $\mu_{1} > 0$ and
$\Gamma_{2}(X_{0}) \cup O$ generates the periodic orbit $L_{2}$
for $\mu_{2} > 0$.

The corresponding domain in $U$, which is the intersection of the
stability regions for $L_{1}$ and $L_{2}$, i.e., i.e., the domain
in which the periodic orbits $L_{1}$ and $L_{2}$ are structurally
stable, will be denoted by $U_{0}$. A stable manifold of $L_{i}$
for the system $X \subset U_{0}$ will be denoted by $W^{s}_{i}$
and the unstable one by $W^{u}_{i}$. If the separatrix value
$A_{i}(X_{0}) > 0$, $W^{u}_{i}$ is a cylinder; if $A_{i}(X_{0}) <
0$, $W^{u}_{i}$ is a M\"obius band. Let us note that, in the case
where ${\cal M}$ is an orientable manifold, $W^{s}_{i}$  will
also be a cylinder if $A_{i}(X_{0}) > 0$. Otherwise it will be a
M\"obius band. However, in the  forthcoming analysis  the signs
of the separatrix values will play an important role. Therefore,
it is natural to distinguish the following three main cases
$$
\begin{array}{lcl}
\mbox { Case~A (orientable)}     & A_{1}(X_{0}) > 0, & A_{2}(X_{0}) > 0,\\
\mbox { Case~B (semiorientable)} & A_{1}(X_{0}) > 0, & A_{2}(X_{0}) < 0,\\
\mbox { Case~C (nonorientable)}  & A_{1}(X_{0}) < 0, & A_{2}(X_{0}) < 0.
\end{array}
$$
In each of the above three cases the domain $U_{0}$ also contains two
bifurcation surfaces $B^{1}_{3}$ and $B^{1}_{4}$:
\begin{enumerate}
\item
In Case~A, $B^{1}_{3}$ corresponds to the inclusion
$\Gamma_{1} \subset W^{s}_{2}$ and $B^{1}_{4}$ corresponds to the inclusion
$\Gamma_{2} \subset W^{s}_{1}$;
\item
In Case~B, $B^{1}_{3}$ corresponds to the inclusion
$\Gamma_{1} \subset W^{s}_{1}$ and $B^{1}_{4}$ corresponds to the inclusion
$\Gamma_{2} \subset W^{s}_{1}$;
\item
In Case~C, along with the above-mentioned generated orbits
$L_{1}$ and $L_{2}$, there also arises a saddle periodic orbit
$L_{3}$ which makes one revolution ``along'' $\Gamma_{1}(X_{0})$
and $\Gamma_{2}(X_{0})$, and if both $W^{u}_{i}$ are M\"obius
bands, $i=1,2$, the unstable manifold $W^{u}_{3}$ of the periodic
orbit $L_{3}$ is a cylinder. In this case the inclusions
$\Gamma_{1} \subset W^{s}_{2}$ and $\Gamma_{2} \subset W^{s}_{3}$
correspond to the surfaces $B^{1}_{3}$ and $B^{1}_{4}$,
respectively.
\end{enumerate}

Suppose that $B^{1}_{3}$ and $B^{1}_{4}$ intersect transversely
over the bifurcational set $B^{2}_{2}$, see Fig.~2.
\begin{figure}
\begin{center}
\includegraphics[width=3 in,height=2.5 in]{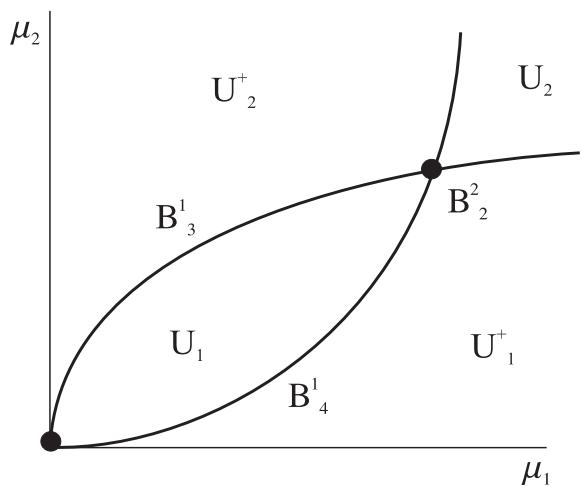}
\end{center}
\caption{$(\mu_{1},\mu_{2})$-bifurcation diagram} \label{fig2}
\end{figure}
 In a
two-parameter family $X(\mu)$ this means that the curves
$B^{1}_{3}$ and $B^{1}_{4}$  intersect at some point $\mu^{1} =
(\mu_{11}, \mu_{12})$. Let us denote a domain lying between
$B^{1}_{3}$ and $B^{1}_{4}$ by $U_{1}$.  Suppose also that for
each $X \in U$ there exists a transversal $D$ with the following
properties:
\begin{enumerate}
\item
The Euclidean coordinates (x, y) can be introduced on $D$ such that
$$
D = \Bigl \{ (x, y): |x| \le 1, |y| < 2 \Bigr \}.
$$
\item
The equation $y=0$ describes a connection component $S$ of the intersection
$W^{s}_{O} \cap D$ such that no $\omega$-semitrajectory that begins on $S$
possesses any point of intersection with $D$ for $t > 0$.
\item
The mapping $T_{1}(X): D_{1} \mapsto D$ and $T_{1}(X): D_{2} \mapsto D$ are
defined along the trajectories of the system $X$, where
$$
D_{1} = \Bigl \{ (x, y): |x| \le 1, 0 < y  \le 1 \Bigr \},
$$
$$
D_{2} = \Bigl \{ (x, y): |x| \le 1, -1 \ge y < 1 \Bigr \},
$$
and $T_{i}(X)$ is written in the form
\begin{equation}\label{p7.4}
\begin{array}{l}
\bar x = f_{i}(x, y), \\
\bar y = g_{i}(x, y),
\end{array}
\end{equation}
where $f_{i}, g_{i} \in {\Bbb C}^{r}$, $i=1,2$.
\item
$f_{i}$ and $g_{i}$ admit continuous extensions on $S$, and
$$
\lim\limits_{y \to 0} f_{i}(x,y) = x^{**}_{i}, \quad
\lim\limits_{y \to 0} g_{i}(x,y) = y^{**}_{i}, \qquad i=1,2.
$$
\item
$$\begin{array}{l}
T_{1} D_{1} \in Pi_{1} = \Bigl \{ (x, y): 1/2 \le x \le 1, |y| <2 \Bigr \},
\\\\
T_{2} D_{2} \in Pi_{2} = \Bigl \{ (x, y): -1 \le x \le -1/2, |y| <2 \Bigr \}.
\end{array}
$$
Let $T(X) \equiv T_{i}(X) \mid D_{i}$, $(f,g) =\equiv (f_{i}, g_{i})$ on
$D_{i}$, $i=1,2$.
\item
Let us impose the following restrictions on $T(X)$

\begin{eqnarray}\label{p7.5}
\left. \begin{array}{ll}
\mbox{(a)} \quad
\| (f_{x}) \| < 1, \\\\
\mbox{(b)} \quad 1 - \| (g_y)^{-1} \| \cdot  \| f_x \| > 2 \sqrt
{ \| (g_y)^{-1} \| \cdot  \| (g_x) \| \cdot \| (g_y)^{-1} \cdot
f_y \| } ,
\\\\
\mbox{(c)} \quad
\|(g_y)^{-1} \| < 1 , \\\\
\mbox{(d)} \quad
\| (g_y)^{-1} \cdot f_y \| \cdot \| g_x \| <
( 1 - \| f_x \|)  (1 - \| (g_y)^{-1} \|).
\end{array} \right \}
\end{eqnarray}
Hereafter, $ \| \cdot \| = \sup\limits_{ (x,y) \in D \backslash S}  \:  | \cdot | $. \\
\end{enumerate}

It follows from the analysis of the behavior of trajectories near
$W^{s}_{O}$ that in a small neighborhood of $S$ the following representation
is valid:
\begin{eqnarray}\label{p7.6}
\begin{array}{l}
f_{1} = x^{**}_{1} + \varphi_{1}(x, y)~~y^{\alpha}, \\
f_{2} = x^{**}_{2} + \varphi_{2}(x, y)(-y)^{\alpha},
\end{array} \qquad
\begin{array}{l}
g_{1} = y^{**}_{1} + \psi_{1}(x, y)~~y^{\alpha}, \\
g_{2} = y^{**}_{2} + \psi_{2}(x, y)(-y)^{\alpha},
\end{array}
\end{eqnarray}
where $\varphi_{1}, \ldots, \psi_{2}$ are smooth with respect to
$x, y$ for $y \neq 0$, and $T_{i}(x)$ satisfies
estimates~(\ref{p7.5}) for sufficiently small $y$. Moreover, the
limit of $\varphi_{1}$ will be denoted by $A_{1}(X)$ and that of
$\psi_{2}$ by $A_{2}(X)$. The functionals $A_{1}(X)$ and
$A_{2}(X)$ will be also called the separatrix values in analogy
with $A_{1}(X_{0})$ and $A_{2}(X_{0})$ which were introduced
above. Let us note that for a system lying in a small
neighborhood of the system $X$ all the conditions~1-6 are
satisfied near $S$. Moreover, the concept of orientable,
semiorientable and nonorientable cases can be extended to any
system $X \in U$. It is convenient to assume, for simplicity, that
$A_{1,2}(X)$ do not vanish. It should be also noted that the point
$P_{i}$ with the coordinates $(x^{**}_{i}, y^{**}_{i})$  is the
first point of intersection of $\Gamma_{i}(X)$ with $D$.

Let us consider the constant
\begin{equation}\label{p7.7}
\begin{array}{ll}
q = \\\\
\displaystyle { 1 + \| f_x \|  \| (g_y)^{-1} \| + \sqrt {1- \|
(g_y)^{-1} \|^{2}  \| (f_x) \| - 4 \| (g_y)^{-1} \|  \| g_x \|
\|(g_{y})^{-1}  f_{y}} \over 2 \|(g_{y})^{-1}\| }.
\end{array}
\end{equation}
Conditions~(\ref{p7.7}) implies that $q > 1$ and hence all the periodic points
will be of the saddle type.

Let $\Sigma$ denote the closure of the set of points of all the trajectories of
the mapping $T(X)$, which are contained entirely in $D$.
$\Sigma$ is described most simply in the domain $U_{1}$. Here the
following theorems hold.

\begin{theorem}\label{t7.1}
If $X \in U_{1}$, $T(X) \mid \Sigma$ is topologically conjugated with the
Bernoulli scheme $(\sigma, \Omega_{2})$ with two symbols.
\end{theorem}

\begin{theorem}\label{t7.2}
The system $X \in U_{2}$, has a two-dimensional limiting set $\Omega$, which
satisfies the following conditions:
\begin{enumerate}
\item
$\Omega$ is structurally unstable.
\item
$[\Gamma_{1} \cup \Gamma_{2} \cap O] \subset \Omega$.
\item
Structurally stable periodic orbits are everywhere dense in $\Omega$.
\item
Under perturbations of $X$ periodic orbits in $\Omega$ disappear as a
result of matching to the saddle separatrix loops $\overline{\Gamma_{1}}$ and
$\overline{\Gamma_{2}}$.
\end{enumerate}
\end{theorem}

Note that in this case the basic periodic orbits will not belong to
$\Omega$. In terms of mappings, the properties of $\Omega$ can be formulated in
more detail. Let us first single out a domain $\tilde D$ on $D$ as follows:
we assume that in Case~A

$$
\tilde D = \Bigr \{ (x,y) \in D_{1} \cup D_{2} \mid y_{2}(x) < y < y_{1}(x)
\Bigr \};
$$
and in Case~B

$$
\tilde D = \Bigr \{ (x,y) \in D_{1} \cup D_{2} \mid y_{12}(x) < y < y_{1}(x)
\Bigr \};
$$
where $y = y_{12}(x)$, $|x| \le 1$, denotes a curve in $D$ whose image lies
on the curve $y=y_{1}(x)$; and finally in Case~C

The closure of points of all the trajectories of the mapping $T(X)$, which
are entirely contained in $\tilde D$, we will denote by $\tilde \Sigma$.

\begin{theorem}\label{t7.3}
Let $X \in U_{2}$. Then:
\begin{enumerate}
\item[I.]
$\tilde \Sigma$ is compact, one-dimensional and consists of two connection
components in Cases~A and C, and of a finite number of connection components
in Case~B.
\item[II.]
$\tilde D$ is foliated by a continuous stable foliation $H^{+}$ into leaves,
satisfying the Lipschitz conditions, along which a point is attracted to
$\tilde \Sigma$; inverse images of the discontinuity line $S:\, y=0$ (with
respect to the mapping $T^{k}$, $k=1,2, \ldots$) are everywhere dense in
$\tilde D$.
\item[III.]
There exits a sequence of $T(X)$-invariant null-dimensional sets $\Delta_{k}$,
$k \in {\Bbb Z_{+}}$, such that $T(X) \mid \Delta_{k}$ is topologically
conjugated with a finite topological Markov chain with a nonzero entropy, the
condition $\Delta_{k} \in \Delta_{k+1}$ being satisfied, and
$\Delta_{k} \to \tilde \Sigma$ as $k \to \infty$.
\item[IV.]
The non-wandering set $\Sigma_{1} \in \tilde \Sigma$ is a closure of
saddle periodic
points of $T(X)$ and either $\Sigma_{1} =\tilde \Sigma$ or
$\Sigma_{1} = \Sigma^{+} \cup \Sigma^{-}$, where:
\end{enumerate}
\begin{enumerate}
\item
$\Sigma^{-}$ is null-dimensional and is an image of the space $\Omega^{-}$ of
a  certain TMC
$(G^{-}, \Omega^{-}, \sigma)$ under the homeomorphism
$\beta: \Sigma^{-} \mapsto \Sigma^{-}$ which conjugates
$\sigma \mid \Omega^{-}$ and  $T(X) \mid \Sigma^{-}$;
$$
\Sigma^{-} = \bigcup\limits_{m=1}^{l(X)} \Sigma^{-}_{m}, l(X) < \infty,
$$
where
$$
T(X) \Sigma^{-}_{m} = \Sigma^{-}_{m}, \qquad
\Sigma^{-}_{m_{1}} \cap \Sigma^{-}_{m_{2}} = \emptyset
$$
for $m_{1} \neq m_{2}$ and $T(X) \mid \Sigma^{-}_{m}$ is transitive;
\item
$\Sigma^{+}$ is compact, one-dimensional and
\item
if $\Sigma^{+} \cap \Sigma^{-} = \emptyset$, $\Sigma^{+}$ is an attracting set in
a certain neighbourhood;
\item
if $\Sigma^{+} \cap \Sigma^{-} \neq \emptyset$, then
$\Sigma^{+} \cap \Sigma^{-} = \Sigma^{+}_{m} \cap \Sigma^{-}_{m}$ for a
certain $m$, and this intersection consists of periodic points of no more
than two periodic orbits, and
\end{enumerate}
\begin{enumerate}
\item[(a)]
if $\Sigma_{m}^{-}$ is finite, $\Sigma^{+}$ is $\omega$-limiting for all the
trajectories in a certain neighbourhood;
\item[(b)]
if $\Sigma_{m}^{-}$ is infinite, $\Sigma^{+}$ is not locally maximal, but
is $\omega$-limiting for all the trajectories in $\tilde D$, excluding those
asymptotic to $\Sigma^{-} \backslash \Sigma^{+}$.
\end{enumerate}
\end{theorem}

Below we will give the conditions under which the existence of the Lorenz
attractor is guaranteed.

Consider a finite-number parameter family of vector field defined by the
system of differential equations
\begin{equation}\label{p7.8}
\dot x = X( x, \mu),
\end{equation}
where $x \in {\Bbb R}^{n+1}$, $\mu \in {\Bbb R}^{m}$, and $X(x, \mu)$
is a ${\Bbb C}^{r}$-smooth functions of $x$ and $\mu$.
Assume that following two conditions hold
\begin{enumerate}
\item[A.]
System~(\ref{p7.8}) has a equilibrium state $O(0,0)$ of the saddle type.
The eigenvalues of the Jacobian at $O(0,0)$ satisfy
$$
\R \lambda_{n} <  \cdots \R \lambda_{2} < \lambda_{1} < 0 < \lambda_{0}.
$$
\item[B.]
The separatrices $\Gamma_{}1$ and $\Gamma_{2}$ of the saddle $O(0,0)$
returns to the origin as $t \to +\infty$.
\end{enumerate}

Then, for  $\mu > 0$ in the parameter space there exists an open
set $V$, whose boundary contains the origin, such that in $V$
system~(\ref{p7.8}) possesses the Lorenz attractor in the following three
cases \cite{SH81}:\\
{\bf Case~1.}
\begin{enumerate}
\item[A]
$\Gamma_{1}$ and $\Gamma$ return to the origin tangentially to each other
along the dominant direction corresponding to the eigenvalue $\lambda_{1}$;
\item[B]
$$
\displaystyle
{1 \over 2} < \gamma < 1, \quad \nu_{i} > 1, \quad
\gamma = -{\lambda_{1} \over \lambda_{0} }, \quad
\nu_{i} = -{ \R \lambda_{i} \over \lambda_{0} };
$$
\item[A]
The separatrix values $A_{1}$ and $A_{2}$ (see above) are equal to zero.
\end{enumerate}
In the general case, the dimension of the parameter space is four since
we may choose $\mu_{1,2}$, to control the behaviour of the separatrices
$\Gamma_{1,2}$ and $\mu_{3,4} = A_{3,4}$. In the case of
the Lorenz symmetry, we need two parameters only.\\
{\bf Case~2.}
\begin{enumerate}
\item[A]
$\Gamma_{1}$ and $\Gamma$ belong to the non-leading manifold
$W^{ss} \in W^{s}$ and enter the saddle along the eigen-direction corresponding
to the real eigenvector $\lambda_{2}$
\item[B]
$$
\displaystyle
{1 \over 2} < \gamma < 1, \quad \nu_{i} > 1, \quad
\gamma = -{\lambda_{1} \over \lambda_{0} }, \quad
\nu_{i} = -{ \R \lambda_{i} \over \lambda_{0} };
$$
\end{enumerate}
In the general case, the dimension of the phase space is equal to four.
Here, $\mu_{3,4}$  control the distance between the separatrices.\\
{\bf Case~3.}
\begin{enumerate}
\item[A]
$\Gamma_{1,2} \notin W^{ss}$;
\item[B]
$\gamma = 1$;
\item[C]
$ A_{1,2} \neq 0$, and $|A_{1,2}|< 2$.
\end{enumerate}
In this case $m=3$, $\mu_{3} = \gamma -1.$

In the case where the system is symmetric, all of these bifurcations are of
codimension~2.

In A.~Shilnikov \cite{AN89,AN93} it was shown that both subclasses~(A) and (C)
are realized in the Shimizu-Morioka model in which the appearance of the Lorenz
attractor and its disappearance through bifurcations of lacunae are explained.
Some systems of type~(A) were studied by Rychlik~\cite{RY90} and those of
type~(C) by Robinson~\cite{RO89}.

The distinguishing features of strange attractors of the Lorenz type is that they
have a complete topological invariant. Geometrically, we can state that
two Lorenz-like attractors are topologically equivalent if the unstable
manifolds of both saddles behave similarly. The formalization of ``similarity''
may be given in terms of {\it kneading invariants} which were introduced by
Milnor and Thurston \cite{MT77} while studying continuous, monotonic mappings
on an interval.
This approach may be applied to certain discontinuous mappings as well. Since
there is a foliation (see above) we may reduce the Poincar\'e map to the
form
\begin{equation}\label{p7.9}
\begin{array}{l}
\bar x = F(x, y), \\
\bar y  = G(y),
\end{array}
\end{equation}
where the right-hand side is, in general, continuous, apart from the
discontinuity line $y=0$, and $G$ is piece-wise monotonic. Therefore, it is
natural to reduce (\ref{p7.9}) to a one-dimensional map
$$
\bar y =G(y),
$$
by using the technique of taking the inverse spectrum,
Guckenheimer and Williams \cite{GW79} showed that a pair of the
kneading invariants is a complete topological invariant for the
associated two-dimensional maps provided $\inf | G^{\prime}| >1$.

\vspace*{-2mm}

\section*{4. Wild strange attractor}

\vskip-5mm \hspace{5mm}

In this section, following the paper by Shilnikov and Turaev \cite{WILD}, we will distinguish a class of dynamical
systems with strange attractors of a new type. The peculiarity of such an attractor is that it may contain a wild
hyperbolic set. We remark that such an attractor is to be understood as an almost stable, chain-transitive closed
set.

Let $X$ be a smooth ($C^r$, $r\geq 4$) flow in $R^n$ ($n\geq 4$)
having an equilibrium state $O$ of {\em a saddle-focus} type with
characteristic exponents $\gamma,-\lambda\pm i\omega,
-\alpha_1,\cdots,-\alpha_{n-3}$ where $\gamma>0$, $0<\lambda< \R
\alpha_j$, $\omega\neq 0$. Suppose
\begin{equation}
\gamma> 2\lambda .\label{cve}
\end{equation}
This condition was introduced in \cite{OSH87} where it was shown,
in particular, that it is necessary in order that no stable periodic orbit
could appear when one of the separatrices of $O$ returns to $O$ as
$t\rightarrow +\infty$ (i.e., when there is {\em a homoclinic loop}; see also
\cite{OSH91}).

Let us introduce coordinates $(x,y,z)$ ($x\in R^1$, $y\in R^2$, $z\in R^{n-3}$)
such that the equilibrium state is in the origin, the one-dimensional unstable
manifold of $O$ is tangent to the $x$-axis and the $(n-1)$-dimensional stable
manifold is tangent to $\{x=0\}$. We also suppose that the coordinates
$y_{1,2}$ correspond to the leading exponents $\lambda\pm i\omega$ and the
coordinates $z$ correspond to the non-leading exponents $\alpha$.

Suppose that the flow possesses {\em a cross-section}, say, the surface
$\Pi:\;\{ \parallel y\parallel=1, \parallel z\parallel\leq1,\mid x\mid\leq1\}$.
The stable manifold $W^s$ is tangent to $\{x=0\}$ at $O$, therefore it is
locally given by an equation of the form $x=h^s(y,z)$ where $h^s$ is a smooth
function $h^s(0,0)=0,\; (h^s)'(0,0)=0$. We assume that it can be written in
such form at least when
$(\parallel y\parallel\leq 1, \parallel z\parallel\leq1)$ and that
$\mid h^s\mid <1$ here. Thus, the surface $\Pi$ is a cross-section for
$W^s_{\rm loc}$ and the intersection of $W^s_{\rm loc}$ with $\Pi$ has the form
$\Pi_0: x=h_0(\varphi,z)$ where $\varphi$ is the angular coordinate:
$y_1=\parallel y\parallel \cos \varphi$, $y_2=\parallel y\parallel\sin\varphi$,
and $h_0$ is a smooth function $-1<h_0<1$. One can make $h_0\equiv 0$ by a
coordinate transformation and we assume that it is done.

We suppose that all the orbits starting on $\Pi\backslash \Pi_0$ return to
$\Pi$, thereby defining the Poincar\'e map:
$T_+:\Pi_+\rightarrow \Pi$ ¨ $T_-:\Pi_-\rightarrow \Pi$,
where $\Pi_+=\Pi\cap\{x>0\}$ ¨ $\Pi_-=\Pi\cap\{x<0\}$. It is evident that
if $P$ is a point on $\Pi$ with coordinates $(x,\varphi,z)$, then
$$
\lim_{x\rightarrow -0}T_-(P)=P^1_-,\,
\lim_{x\rightarrow +0}T_+(P)=P^1_+,$$
where $P^1-$ and $P^1_+$ are the first intersection points of the
one-dimensional separatrices of $O$ with $\Pi$. We may therefore define the
maps $T_+$ and $T_-$ so that
\begin{equation}
T_-(\Pi_0)=P^1_-,\,T_+(\Pi_0)=P^1_+ .\label{cc}
\end{equation}

Evidently, the region ${\cal D}$ filled by the orbits starting on
$\Pi$ (plus the point $O$ and its separatrices) is {\em an
absorbing domain} for the system $X$ in the sense that the orbits
starting in $\partial {\cal D}$ enter ${\cal D}$ and stay there
for all positive values of time $t$. By construction, the region
${\cal D}$ is the cylinder $\{\parallel y\parallel\leq1,\parallel
z\parallel\leq1,\mid x\mid\leq1\}$ with two glued handles
surrounding the separatrices, see Fig.3.
\begin{figure}
\begin{center}
\includegraphics[width=3 in,height=2.5 in]{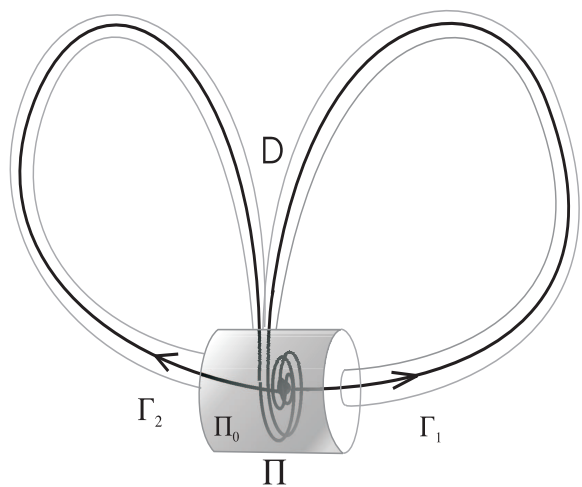}
\end{center}
\caption{Construction of the wild attractor} \label{fig3}
\end{figure}

We suppose that the (semi)flow is {\em pseudohyperbolic} in ${\cal D}$.
It is convenient for us to give this notion a sense more strong than it is
usually done \cite{HPS77}. Namely, we propose the following\\
{\bf Definition.} {\em A semi-flow is called {\em pseudohyperbolic} if the
following two conditions hold:\\
\begin{enumerate}
\item[{\bf A}]
At each point of the phase space, the tangent space is uniquely
decomposed (and this decomposition is invariant with respect to the linearized
semi-flow) into a direct sum of two subspaces $N_1$ and $N_2$ (continuously
depending on the point) such that the maximal Lyapunov exponent in $N_1$
is strictly less than the minimal Lyapunov exponent in $N_2$: at each point
$M$, for any vectors $u\in N_1(M)$ and $v\in N_2(M)$
$$\limsup_{t\rightarrow+\infty}\frac{1}{t}\ln\frac{\|u_t\|}{\|u\|}<
\liminf_{t\rightarrow+\infty}\frac{1}{t}\ln\frac{\|v_t\|}{\|v\|}$$
where $u_t$ and $v_t$ denote the shift of the vectors $u$ and $v$ by the
semi-flow linearized along the orbit of the point $M$;\\
\item[{\bf B}]
The linearized flow restricted on $N_2$ is volume expanding:
$$V_t \ge const\cdot e^{\sigma t} V_0$$
with some $\sigma>0$; here, $V_0$ is the volume of any region in $N_2$ and
$V_t$ is the volume of the shift of this region by the linearized semi-flow.
\end{enumerate} }

The additional condition {\bf B} is new here and it prevents of appearance of
stable periodic orbits. Generally, our definition includes the case where the
maximal Lyapunov exponent in $N_1$ is non-negative everywhere. In that case,
according to condition {\bf A}, the
linearized semi-flow is expanding in $N_2$ and
condition {\bf B} is satisfied trivially. In the present paper we consider
the opposite case where the linearized semi-flow is {\em exponentially
contracting} in $N_1$, so condition {\bf B} is essential here.

Note that the property of pseudo-hyperbolicity is stable with respect to small
smooth perturbation of the system: according to \cite{HPS77} the invariant
decomposition of the tangent space is not destroyed by small perturbations
and the spaces $N_1$ and $N_2$ depend continuously on the system. Hereat,
the property of volume expansion in $N_2$ is also stable with respect to small
perturbations.

Our definition is quite broad; it embraces, in particular, hyperbolic flows
for which one may assume $(N_1,N_2)=(N^s,N^u\oplus N_0)$ or
$(N_1,N_2)=(N^s\oplus N_0,N^u)$ where $N^s$ and $N^u$ are, respectively,
the stable and unstable invariant subspaces and $N_0$ is a one-dimensional
invariant subspace spanned by the phase velocity vector. The geometrical
Lorenz model from \cite{ABSH77,ABSH83} or \cite{GW79} belongs also to this
class: here $N_1$ is tangent to the contracting invariant foliation of
codimension two and the expansion of areas in a two-dimensional subspace $N_2$
is provided by the property that
the Poincar\'e map is expanding in a direction transverse to the contracting
foliation.

In the present paper we assume that $N_1$ has codimension three:
$dim N_1=n-3$ and $dim N_2=3$ and that the linearized flow (at $t\geq 0$) is
exponentially contracting on $N_1$. Condition {\bf A} means here that if for
vectors of $N_2$ there is a contraction, it has to be weaker than those on
$N_1$. To stress the last statement, we will call $N_1$ {\em the strong stable
subspace} and $N_2$ {\em the center subspace} and will denote them as
$N^{ss}$ and $N^c$ respectively.

We also assume that the coordinates $(x,y,z)$ in $R^n$ are such that
at each point of ${\cal D}$ the space $N^{ss}$ has a non-zero projection
onto the coordinate space $z$, and $N^{c}$ has a non-zero projection onto the
coordinate space $(x,y)$.

Note that our pseudohyperbolicity conditions are satisfied at the point $O$
from the very beginning:
the space $N^{ss}$ coincides here with the coordinate space $z$, and $N^{c}$
coincides with the space $(x,y)$; it is condition (\ref{cve}) which
guarantees the expansion of volumes in the invariant subspace $(x,y)$.
The pseudohyperbolicity of the linearized flow is automatically inherited
by the orbits in a small neighborhood of $O$. Actually, we require that
this property would extend into the non-small neighborhood ${\cal D}$ of $O$.

According to \cite{HPS77}, the exponential contraction in $N^{ss}$
implies the existence of an invariant contracting foliation ${\cal
N}^{ss}$ with $C^r$-smooth leaves which are tangent to $N^{ss}$.
As in \cite{ABSH83}, one can show that the foliation is absolutely
continuous. After a factorization along the leaves, the region
${\cal D}$ becomes a branched manifold (since ${\cal D}$ is
bounded and the quotient-semiflow expands volumes it follows
evidently that the orbits of the quotient-semiflow must be glued
on some surfaces in order to be bounded; cf.\cite{W79}).

The property of pseudohyperbolicity is naturally inherited by the Poincar\'e
map $T\equiv (T_+,T_-)$ on the cross-section $\Pi$: here, we have: \\
\begin{enumerate}
\item[{\bf A}$^*$]
There exists a foliation with smooth leaves of the form
$(x,\varphi)=h(z)\mid_{-1\leq z\leq 1}$, where the derivative $h'(z)$ is
uniformly bounded, which possesses the following properties: the foliation is
{\em invariant} in the sense that if $l$ is a leaf, then
$T_+^{-1}(l\cap T_+(\Pi_+\cup\Pi_0))$ and $T_-^{-1}(l\cap T_-(\Pi_-\cup\Pi_0))$
are also leaves of the foliation (if they are not empty sets); the
foliation is {\em absolutely continuous} in the sense that the projection
along the leaves from one two-dimensional transversal to another increases or
decreases the area in a finite number of times and the coefficients of
expansion or contraction of areas are bounded away from zero and infinity;
the foliation is {\em contracting} in the sense that if two points belong
to one leaf, then the distance between the iterations of the points with
the map $T$ tends to zero exponentially;
\item[{\bf B}$^*$]
The quotient maps $\tilde{T_+}$ and $\tilde{T_-}$ are
area-expanding.
\end{enumerate}

\begin{statement}  Let us write the map $T$ as
$$(\bar{x},\bar{\varphi})=g(x,\varphi,z),\;\; \bar{z}=f(x,\varphi,z),$$
where $f$ and $g$ are functions smooth at $x\neq 0$ and discontinuous at $x=0$:
$$\lim_{x\rightarrow -0}(g,f)=(x_-,\varphi_-,z_-)\equiv P^1_-,
\;\;\lim_{x\rightarrow +0}(g,f)=(x_+,\varphi_+,z_+)\equiv P^1_+\;.$$
Let
\begin{equation}
\det \;\frac{\partial g}{\partial (x,\varphi)}\neq 0\;.
\label{ct0}
\end{equation}
Denote
$$\begin{array}{ll}
A=\frac{\partial f}{\partial z}-\frac{\partial f}{\partial (x,\varphi)}
\left(\frac{\partial g}{\partial (x,\varphi)}\right)^{-1}
\frac{\partial g}{\partial z}, &
B=\frac{\partial f}{\partial (x,\varphi)}
\left(\frac{\partial g}{\partial (x,\varphi)}\right)^{-1} ,\\
C=\left(\frac{\partial g}{\partial (x,\varphi)}\right)^{-1}
\frac{\partial g}{\partial z}, &
D=\left(\frac{\partial g}{\partial (x,\varphi)}\right)^{-1}.
\end{array}$$
If
\begin{equation}
\lim_{x\rightarrow 0} C=0,\;\; \lim_{x\rightarrow 0} \parallel A
\parallel \; \parallel D \parallel =0 ,\label{ct1}
\end{equation}
\begin{equation}
\sup_{P\in \Pi\backslash \Pi_0}
\sqrt{\parallel A \parallel \; \parallel D \parallel} +
\sqrt{\sup_{P\in \Pi\backslash \Pi_0}\parallel B \parallel
\sup_{P\in \Pi\backslash \Pi_0}\parallel C \parallel} <1,
\label{ct2}
\end{equation}
then the map has a continuous invariant foliation with smooth leaves of the
form $(x,\varphi)=h(z)\mid_{-1\leq z\leq 1}$ where the derivative $h'(z)$ is
uniformly bounded. If, additionally,
\begin{equation}
\sup_{P\in \Pi\backslash \Pi_0}\parallel A \parallel +
\sqrt{\sup_{P\in \Pi\backslash \Pi_0}\parallel B \parallel
\sup_{P\in \Pi\backslash \Pi_0}\parallel C \parallel} <1,
\label{ct3}
\end{equation}
then the foliation is contracting and if, moreover, for some
$\beta>0$
\begin{equation}
\begin{array}{l}
\mbox{the functions } A\mid x\mid^{-\beta},\;D\mid x\mid^{\beta},\; B,\;C
\mbox{ are uniformly bounded and} \\ \mbox {H\"older continuous,}\\
\mbox{and }\displaystyle\frac{\partial \ln \det\; D}{\partial z}
\mbox{ and    }
\frac{\partial \ln \det\; D}{\partial (x,\varphi)}D\mid x\mid^\beta
\mbox{ are uniformly bounded,}
\end{array}
\label{ct4}
\end{equation}
then the foliation is absolutely continuous. The additional condition
\begin{equation}
\sup_{P\in \Pi\backslash \Pi_0}\sqrt{\det\; D} +
\sqrt{\sup_{P\in \Pi\backslash \Pi_0}\parallel B \parallel
\sup_{P\in \Pi\backslash \Pi_0}\parallel C \parallel} <1
\label{ct5}
\end{equation}
guarantees that the quotient map $\tilde{T}$ expands areas.
\end{statement}

It follows from \cite{OSH87,OSH91} that in the case where the
equilibrium state is a saddle-focus, the Poincar\'e map near
$\Pi_0=\Pi\cap W^s$ is written
in the following form under some appropriate choice of the coordinates.
\begin{equation}
(\bar{x},\bar{\varphi})=Q_\pm(Y,Z),\;\;\; \bar{z}=R_\pm(Y,Z).
\label{glm}
\end{equation}
Here
\begin{equation}
\begin{array}{l}
Y=\mid x \mid^\rho \left(\begin{array}{ll}
\cos (\Omega\ln \mid x\mid +\varphi)& \sin (\Omega\ln \mid x\mid+\varphi)\\
-\sin (\Omega\ln \mid x\mid +\varphi)& \cos (\Omega\ln \mid x\mid+\varphi)
\end{array}\right)+\Psi_1(x,\varphi,z),\\
\;\\
Z=\Psi_2(x,\varphi,z),
\end{array}
\label{lcm}
\end{equation}
where $\rho=\lambda/\gamma<1/2$ (see ({\ref{cve})), $\Omega=\omega/\gamma$ and,
for some $\eta>\rho$,
\begin{equation}
\parallel\frac{\partial^{p+\mid q \mid} \Psi_i}
{\partial x^p \partial (\varphi,z)^q}\parallel = O(\mid x \mid^{\eta - p}),\;\;
0\leq p+\mid q\mid \leq r-2\; ;
\label{oc}
\end{equation}
the functions $Q_\pm$, $R_\pm$ in (\ref{glm}) (``+'' corresponds to $x>0$ -
the map $T_+$, ``-'' corresponds to $x<0$ - the map $T_-$) are smooth
functions in a neighborhood of $(Y,Z)=0$ for which the Taylor expansion can
be written down
\begin{equation}
Q_\pm=(x_\pm,\varphi_\pm)+a_\pm Y + b_\pm Z + \cdots\;,\;\;\;
R_\pm=z_\pm+c_\pm Y + d_\pm Z + \cdots\ . \label{gle}
\end{equation}

It is seen from (\ref{glm})-(\ref{gle}) that if $O$ is a saddle-focus
satisfying (\ref{cve}), then if $a_+\neq 0$ and $a_-\neq 0$, the map $T$
satisfies conditions (\ref{ct1}) and (\ref{ct4}) with $\beta \in (\rho,\eta)$.
Furthermore, analogues of conditions (\ref{ct0}),(\ref{ct2}),(\ref{ct3}),
(\ref{ct5}) are fulfilled where the supremum should be taken
not over $\mid x \mid \leq 1$ but it is taken over small $x$. The map
(\ref{glm}),(\ref{lcm}),(\ref{gle}) is easily continued onto the whole
cross-section $\Pi$ so that the conditions of the lemma were fulfilled
completely. An example is given by the map
\begin{equation}
\begin{array}{l}
\bar{x}=0.9 \mid x \mid^\rho \cos (\ln \mid x\mid +\varphi),\\
\;\\
\bar{\varphi}=3 \mid x \mid^\rho \sin (\ln \mid x\mid +\varphi),\\
\;\\
\bar{z}=(0.5 + 0.1 z \mid x \mid^\eta)\;{\rm sign}\; x
\end{array}
\end{equation}
where $0.4=\rho<\eta$

As stated above, the expansion of volumes by the quotient-semiflow
restricts the possible types of limit behavior of orbits. Thus, for instance,
{\em in $\cal D$ there may be no stable periodic orbits}. Moreover, {\em
any orbit in $\cal D$ has a positive maximal Lyapunov exponent}. Therefore,
one must speak about a strange attractor in this case.

Beforehand, we recall some definitions and simple facts from
topological dynamics. Let $X_tP$ be the time-$t$ shift of a point
$P$ by the flow $X$. For given  $\varepsilon>0$ and $\tau>0$ let
us define as {\em an $(\varepsilon,\tau)$-orbit} as a sequence of
points $P_1,P_2,\cdots,P_k$ such that $P_{i+1}$ is at a distance
less than $\varepsilon$ from $X_tP_i$ for some $t>\tau$. A point
$Q$ will be called $(\varepsilon,\tau)${\em-attainable} from $P$
if there exists an $(\varepsilon,\tau)$-orbit connecting $P$ and
$Q$; and it will be called {\em attainable} from $P$ if, for some
$\tau>0$, it is $(\varepsilon,\tau)$-attainable from $P$ for any
$\varepsilon$ (this definition, obviously, does not depend on the
choice of $\tau>0$). A set $C$ is attainable from $P$ if it
contains a point attainable from $P$. A point $P$ is called {\em
chain-recurrent} if it is attainable from $X_tP$ for any $t$. A
compact invariant set $C$ is called {\em chain-transitive} if for
any points $P\in C$ and $Q\in CC$ and for any $\varepsilon>0$ and
$\tau>0$ the set $C$ contains an $(\varepsilon,\tau)$-orbit
connecting $P$ and $Q$. Clearly, all points of a chain-transitive
set are chain-recurrent.

A compact invariant set $C$ is called {\em orbitally stable},
if for any its neighborhood $U$ there is a neighborhood $V(C)\subseteq U$
such that the orbits starting in $V$ stay in $U$ for all $t\geq 0$.
An orbitally stable set will be called {\em completely stable} if for any its
neighborhood $U(C)$
there exist $\varepsilon_0>0$, $\tau>0$ and a neighborhood $V(C)\subseteq U$
such that the $(\varepsilon_0,\tau)$-orbits starting in $V$ never leave $U$.
It is known, that a set $C$ is orbitally stable if and only if
$\displaystyle C=\bigcap\limits_{j=1}^\infty U_j$ where $\{U_j\}_{j=1}^\infty$
is a system of embedded open invariant (with respect to the forward flow) sets,
and $C$ is completely stable if the sets $U_j$ are not just invariant but
they are absorbing domains (i.e.; the orbits starting on $\partial U_j$ enter
inside $U_j$ for a time interval not greater than some $\tau_j$; it is clear in this
situation that $(\varepsilon,\tau)$-orbits starting on
$\partial U_j$ lie always inside $U_j$ if $\varepsilon$ is sufficiently small
and $\tau\geq \tau_j$). Since the maximal invariant set ({\em the maximal
attractor}) which lies in any absorbing domain is, evidently,
{\em asymptotically stable}, it follows that any completely stable set
is either asymptotically stable or is an intersection of a countable number of
embedded closed invariant asymptotically stable sets.

\noindent {\bf Definition.} {\em We call the set $\cal A$ of the
points attainable from the equilibrium state $O$ {\em the
attractor} of the system $X$.}

This definition is justified by the following theorem.
\begin{theorem}\label{t9.1}
The set $\cal A$ is chain-transitive, completely stable and attainable
from any point of the absorbing domain $\cal D$.
\end{theorem}

Let us consider a one-parameter family $X_\mu$ of such systems assuming that:\\
{\em a homoclinic loop of the saddle-focus $O$ exists at $\mu=0$,}\\
i.e., one of the separatrices (say, $\Gamma_+$) returns to $O$ as
$t\rightarrow+\infty$. In other words, we assume that the family $X_\mu$
intersects, at $\mu=0$, a bifurcational surface filled by systems with
a homoclinic loop of the saddle-focus and we suppose that this intersection is
{\em transverse}.
The transversality means that when $\mu$ varies, the loop splits and if $M$
is the number of the last point of intersection of the separatrix $\Gamma_+$
with the cross-section $\Pi$ at $\mu=0$ ($P_M^+\in\Pi_0$ at $\mu=0$), then the
distance between the point $P_M^+$ and $\Pi_0$ changes with a
``non-zero velocity'' when $\mu$ varies. We choose the sign of $\mu$
so that $P_M^+\in \Pi_+$ when $\mu>0$ (respectively, $P_M^+\in \Pi_-$
when $\mu<0$).
\begin{theorem}\label{t9.4}
There exists a sequence of intervals $\Delta_i$ (accumulated at $\mu=0$) such
that when $\mu\in\Delta_i$, the attractor ${\cal A}_\mu$ contains a wild
set (non-trivial transitive closed hyperbolic invariant set whose unstable
manifold has points of tangency with its stable manifold). Furthermore,
for any $\mu^*\in\Delta_i$, for any system $C^r$-close to a system $X_\mu^*$,
its attractor $\cal A$ also contains the wild set.
\end{theorem}

We have mentioned earlier that the presence of structurally unstable
(nontransverse) homoclinic trajectories leads to non-trivial
dynamics. Using results \cite{GTS1,GTS93}  we can conclude  that
the systems whose attractors contain structurally non-transverse
homoclinic trajectories as well as structurally unstable periodic
orbits of higher orders of degeneracies are dense in the given regions
in the space of dynamical systems. In particular, the values of
$\mu$ are dense in the intervals $\Delta_{i}$ for which an
attractor of the system contain a periodic orbit of the
saddle-saddle type along with its three-dimensional unstable
manifold. For these parameter values, the topological dimension of
such an attractor is already not less than three. The latter
implies that the given class of systems is an example of
hyperchaos.

\vspace*{-2mm}

\section*{5. Summary} \addsec

\vskip-5mm \hspace{5mm}

The above listed attractors are, in the ideal, the most suitable
images of dynamical chaos. Even though some of them are
structurally unstable, nevertheless it is important that no stable
periodic orbits appear in the system under small smooth
perturbations. Nonetheless, excluding the attractors of the Lorenz
type, no others have been ever observed in nonlinear dynamics so
far. A research has frequently to deal with the models in which
despite the complex behavior of the trajectories appears to be so
visually convincing, nevertheless explicit statements regarding
the exponential instability of the trajectories in the limit set
can be debatable, and therefore should be made with caution. In
numeric experiments with such model one finds a positive Lyapunov
exponent, a continuous frequency spectrum, fast decaying
correlation functions etc., i.e. all the attributes of dynamical
chaos seem to get fulfilled so that the presence of the dynamical
chaos causes no doubts. However, this ``chaotic attractor" may and
often do contain countably many stable periodic orbits which have
long periods and week and narrow attraction basins. Besides the
corresponding stability regions are relatively miniature in the
parameter space under consideration and whence those orbits do
not reveal themselves ordinarily in numeric simulations except
some quite large stability windows where they are clearly visible.
If it is the case, the quasi-attractor \cite{ASH83} is a more
appropriate term for such chaotic set. The natural cause for this
rather complex dynamics is homoclinic tangencies.  Today the the
systems with homoclinic tangencies are the target of many studies.
We briefly outline some valuable facts proven for 3D systems and
2D diffeomorphisms. It will be clear that these results will also
hold for the general case where there may be some other
peculiarities, as for instance the co-existence of countable sets
of saddle periodic orbits of distinct topological types (see
\cite{GTS93}).

We suppose that the system possesses an absorbing area embracing the
hyperbolic basis set in which  the stable and unstable
subsets may touch each other. If it is so, such a hyperbolic set is
called {\em wild}. It follows then that either the system itself
or a close one will have a saddle periodic orbit with
non-transverse homoclinic trajectory along which the stable and
unstable manifolds of the cycle have the tangency. In general, the
tangency is quadratic. Let the saddle value $|\lambda \,
\gamma|$ be less then 1, where $\lambda$ and $\gamma$ are the
multipliers  of the saddle periodic orbit. This condition is always
true when
the divergence of the vector field is negative in the  absorbing
area.  Therefore, near the given system there will exist the
so-called Newhouse regions \cite{N79} in the space of the
dynamical system, i.e. the regions of dense structural
instability. Moreover, a system in the Newhouse region has
countably many stable periodic orbits which cannot principally
be separated from the hyperbolic subset. If additionally
this hyperbolic set contains a saddle periodic orbit with
the saddle value exceeding one, then there will be a countable set of
repelling periodic orbits next to it, and whose closure is not
separable from the hyperbolic set either. The pictures becomes
ever more complex if the divergence of the vector field is
sign-alternating in the absorbing area. Such exotic dynamics
requires infinitely many continuous topological invariant ---
{\it  moduli}, needed for the proper description of the system in the
Newhouse regions. This result comes from the fact that
the systems with the countable set of periodic orbits of arbitrary high
degrees of degeneracies are dense in the Newhouse regions
\cite{GTS92,GTS93}. That is why we ought to conclude in a bitter way:
the complete theoretical analysis of the models, which admit homoclinic
tangencies, including complete bifurcation diagrams and so forth
is non realistic.

\end{document}